\documentclass[11pt,a4paper]{article}
\usepackage[ansinew]{inputenc}
\usepackage{amsmath,amssymb,graphicx,latexsym,amsfonts}
\usepackage{pstricks}
\usepackage{graphicx,rotating}
\usepackage{paralist}
\usepackage{hyperref}

\newtheorem{satz}{Theorem}
\newtheorem{lemma}{Lemma}

\newtheorem{bemerkung}{Remark}
\makeatletter
\renewenvironment{description}{\list{}{\setlength{\leftmargin}{1em}%
 \labelwidth\z@\itemindent-\leftmargin}}{\endlist}

\makeatother

\newcommand{\Hess}{\operatorname{Hess}}
\newcommand{\Ar}{\operatorname{Ar}}

\newcommand{\rank}{\operatorname{rank}}

\newcommand{\rb}[1]{\raisebox{1.5ex}[-1.5ex]{#1}}

\title{A Quintic Hypersurface in $\mathbb{P}^8(\mathbb{C})$ with Many Nodes}
\author{Oliver Schmidt\footnote{Fraunhofer--Institut für Techno- und
    Wirtschaftsmathematik Kaiserslautern. E-Mail: 
\url{oliver.schmidt@itwm.fraunhofer.de}} \and 
Oliver Labs\footnote{Universit\"at des Saarlandes, Saarbr\"ucken.
E-Mail: \url{Labs@math.uni-sb.de}} \and 
Duco van Straten\footnote{Johannes Gutenberg--Universit\"at Mainz.
Partly supported by DFG Sonderforschungsbereich/Transregio 45.
E-Mail: \url{straten@mathematik.uni-mainz.de}}}
\date{} 

\newcommand{\dontshow}[1]{}

\begin{document}

\maketitle
\thispagestyle{empty}


\begin{abstract}
\noindent
We construct a hypersurface of degree $5$ in projective space
$\mathbb{P}^8(\mathbb{C})$ which contains exactly 23436 ordinary
nodes and no further singularities. 
This limits the maximum number $\mu_{8}(5)$ of ordinary nodes a
hyperquintic in $\mathbb{P}^8(\mathbb{C})$ can have to  
 $ 23436  \leq  \mu_{8}(5)  \leq  27876$.
Our method generalizes the approach by the $3^{\text{rd}}$ author
for the construction of a quintic threefold with $130$ nodes
in an earlier paper. 
\end{abstract}

\section*{Introduction}
Let $\mu_{n}(d)$ be the maximum number of ordinary nodes a
hypersurface of degree $d$ in $\mathbb{P}^n:=\mathbb{P}^n(\mathbb{C})$
can have. It is known only for a few nontrivial cases:
For curves in the plane we have $\mu_{2}(d)=d(d-1)/2$.
In three--space, $\mu_3(d)$ is only known for $d\le 6$; see \cite{Bar,JR}
for the case of degree six and \cite{Lab} for an extensive overview.
In $\mathbb{P}^n$ with $n\geq4$, the best known upper bound is Varchenko's
\emph{spectral bound}  \cite{Var} 
\[ \mu_{n}(d) \leq \Ar_n(d), \]
where $\Ar_n(d)$ is Arnold's number:
\[
  \Ar_n(d):=
    \#\Big\{(k_0,\hdots,k_n)\in\big((0,d)\cap\mathbb{Z}\big)^{n+1}
            \,\Big|\,
            \sum_{i=0}^n k_i=\Big\lfloor\frac{nd}{2}\Big\rfloor+1 \Big\}.
\]
All currently best known lower bounds follow from symmetric
constructions: 
Kalker \cite{Kal} constructed $\Sigma_{n}$-symmetric cubics 
which show $\mu_{n}(3)=\Ar_n(3)={n+1\choose\lfloor\frac{n}{2}\rfloor}$
for any $n$. 
Goryunov constructed $A_{n+1}$- and $B_{n+1}$-symmetric quartics
in $\mathbb{P}^n$, which reach approximately $86\%$ of the
Arnold-Varchenko upper bound (cf.~\cite{Gor}).  
In \cite{vStr}, a $\Sigma_6$-symmetric quintic in $\mathbb{P}^4$ with
$130$ nodes was constructed which limits the possibilities for
$\mu_{4}(5)$ to $130 \leq \mu_{4}(5) \leq 135 = \Ar_4(5).$

In sections \ref{sec:Sigmanp2HQ} to \ref{sec:ordinaryNodes}, we
consider the case of $\Sigma_{n+2}$-invariant 
quintics and construct an example in $\mathbb{P}^8$ with  
$23436$ nodes which yields
\[
 23436 \leq \mu_8(5) \leq 27876 = \Ar_8(5).
\]
For most other $n$, it seems that a pentagon--symmetric construction
yields more nodes than our approach; we
discuss this briefly in section \ref{sec:conclRem}. 

\section[$\Sigma_{\lowercase{n}+2}$-symmetric Hyperquintics]
        {$\Sigma_{n+2}$-symmetric Hyperquintics}
  \label{sec:Sigmanp2HQ} 
Adapting the approach used in \cite{vStr}, we consider the 1-pa\-ra\-me\-ter-family of $\Sigma_{n+2}$-symmetric hyperquintics \mbox{$Q:=Q_{(\alpha:\beta)}$} given by
\[
 F_{(\alpha:\beta)}:=\alpha S_5 + \beta S_2 S_3=0\,,\quad (\alpha:\beta)\in\mathbb{P}^1\,,
\]
in projective space $\mathbb{P}^n(\mathbb{C})$, which is defined by $S_1=0$ in $\mathbb{P}^{n+1}(\mathbb{C})$.
Here, $S_i$ denotes the $i$-th elementary-symmetric polynomial in the space coordinates of $\mathbb{P}^{n+1}$:
\[
 S_i=\sum_{0\leq j_1<\hdots<j_i\leq n+1} x_{j_1}\cdot\hdots\cdot x_{j_i}, \qquad i=1,\hdots,5\,.
\]
To determine the sin\-gu\-lar locus of each quintic in the pencil, 
it turns out to be convenient to rewrite $F_{(\alpha:\beta)}$ in terms of the \emph{$i$-th power sums} in the coordinates $x_j$ defined by 
\[
 C_i := \sum^{n+1}_{j=0} x_j^i,\qquad i=1,\hdots,5.
\]

Modulo $S_1$, we have the following identities:
\begin{align*}
 S_1 & = C_1,\\
 S_2 & = \textstyle -\frac{1}{2} C_2,\\
 S_3 & = \textstyle \frac{1}{3} C_3,\\
 S_4 & = \textstyle -\frac{1}{4} C_4 + \frac{1}{8} C_2^2,\\
 S_5 & = \textstyle \frac{1}{5} C_5 - \frac{1}{6} C_2 C_3.
\end{align*}
So the hyperquintic $Q=Q_{(\alpha:\beta)}$ is given by 
\[
 F_{(\alpha:\beta)} =
    \alpha S_5+\beta S_2 S_3=\textstyle\frac{\alpha}{5} C_5-\frac{\alpha+\beta}{6} C_2 C_3=0\,.
\]
Since $F_{(0:1)} = -\frac{1}{6} C_2 C_3 = S_2 S_3$ clearly has the 
projective variety $S_2=S_3=0$ as singular locus, we assume
$\alpha\not=0$. 
The singular points of the hyperquintics are those where the gradients of the defining equations in $\mathbb{P}^{n+1}$ are dependent. So we have
\begin{align*}
 \eta \text{ singular} \qquad
   & \Leftrightarrow   \qquad
     \rank\left(\begin{array}{ccc}
       \partial_0 F_{(\alpha:\beta)}(\eta) & \hdots & \partial_{n+1} F_{(\alpha:\beta)}(\eta)\\
       \partial_0 S_1(\eta)                & \hdots & \partial_{n+1} S_1(\eta)\\
                \end{array}\right) \leq1\\
   & \Leftrightarrow   \qquad
     \rank\left(\begin{array}{ccc}
       \partial_0 F_{(\alpha:\beta)}(\eta) & \hdots & \partial_{n+1} F_{(\alpha:\beta)}(\eta)\\
                 1                         & \hdots & 1\\
                \end{array}\right) \leq1\\
   & \Leftrightarrow   \qquad
                     \exists \,\mu\in\mathbb{C}:
                     \partial_i F_{(\alpha:\beta)}(\eta) = \mu\,, \qquad i=0,\hdots,n+1\,.
\end{align*}

Hence, for all indices $i=0,\hdots,n+1$ we obtain 
\[
 \sum_{j=0}^{n+1}\partial_j F_{(\alpha:\beta)}(\eta)=(n+2)\cdot\mu
 =(n+2)\cdot\partial_i F_{(\alpha:\beta)}(\eta),
\]
which leads via $S_1=0$ to the following lemma. 
\begin{lemma}
\label{Vorarbeiten}
Each coordinate $\eta_i$ of a singularity $\eta$ of the hyperquintic $Q_{(\alpha:\beta)}$
in \mbox{$\mathbb{P}^n=V(S_1(x_0,\hdots,x_{n+1}))$} 
is a root of 
\[
 P(X):=P_{\lambda}(X):=\textstyle X^4-X^2\cdot\lambda C_2 -X\cdot\frac{2}{3}\lambda C_3
     +\frac{1}{n+2} \big(\lambda C_2^2-C_4\big)=0,
\]
where $\lambda:=\frac{\alpha+\beta}{2\alpha}$.
\end{lemma}
\vspace{-7mm}
\begin{flushright}$\square$\end{flushright}
Note that the sum of the four roots of $P(X)$ is zero since the term
$X^3$ does not occur. 

\section{The Family of $\Sigma_{10}$-symmetric Hyperquintics in $\mathbb{P}^8$}
\label{family_hyperquintics}
We now specialize to the case $n=8$. 
According to Lemma \ref{Vorarbeiten}, each coordinate $\eta_i$ of a singularity $\eta$ of the hyperquintic $Q_{(\alpha:\beta)}$ in $\mathbb{P}^8$ satisfies 
\[
 P(X)=\textstyle X^4-X^2\cdot\lambda C_2 -X\cdot\frac{2}{3}\lambda C_3
     +\frac{1}{10} \big(\lambda C_2^2-C_4\big)=0\,,
\]
where 
$\lambda=\frac{\alpha+\beta}{2\alpha}$. 
A priori, there are 23 cases to check, since the 10 coordinates may be distributed over the four roots $a,b,c,d$ of $P$ as follows: 

\begin{center}
 \begin{tabular}{|ll|ll|ll|}\hline
  Case 1: & $10a          $ & Case 9:  & $6a, 3b,  c    $ & Case 17: & $4a, 4b, 2c     $\\
  Case 2: & $9a,  b       $ & Case 10: & $6a, 2b, 2c    $ & Case 18: & $4a, 4b,  c,  d $\\
  Case 3: & $8a, 2b       $ & Case 11: & $6a, 2b,  c, d $ & Case 19: & $4a, 3b, 3c     $\\
  Case 4: & $8a,  b, c    $ & Case 12: & $5a, 5b        $ & Case 20: & $4a, 3b, 2c,  d $\\
  Case 5: & $7a, 3b       $ & Case 13: & $5a, 4b,  c    $ & Case 21: & $4a, 2b, 2c, 2d $\\
  Case 6: & $7a, 2b, c    $ & Case 14: & $5a, 3b, 2c    $ & Case 22: & $3a, 3b, 3c,  d $\\
  Case 7: & $7a,  b, c, d $ & Case 15: & $5a, 3b,  c, d $ & Case 23: & $3a, 3b, 2c, 2d $\\
  Case 8: & $6a, 4b       $ & Case 16: & $5a, 2b, 2c, d $ & &                  \\ \hline
 \end{tabular}
\end{center}

We analyse some example cases here; the remaining cases can be found in the appendix. 
First, we determine only the $\Sigma_{10}$-orbit length of the corresponding singularity $\eta$. Then, we further check for nodes in those cases that produced the longest orbits under $\Sigma_{10}$.

\begin{description}
 \item[Case 1] does not occur, since on the one hand $\eta=(x:\hdots:x)\in\mathbb{P}^8$, 
      and on the other hand the sum of its coordinates has to be zero. 
 \item[Case 2] Assume that $\eta = (1:1:1:1:1:1:1:1:1:-9)$. Hence $C_2=90$, $C_3=-720$, 
      $C_4=6570$, and 
      \[
       P(X) = \textstyle X^4-90\lambda X^2 +480\lambda X +810\lambda-657\,.
      \]
      Requiring $P(1)=P(-9)=0$, we obtain $\lambda=\frac{41}{75}$, thus 
      \[
       (\alpha:\beta)=(75:7)\,.
      \]
      The length of the $\Sigma_{10}$-orbit of $\eta$ is 10.
 \item[Case 7] A priori, we have $\eta=(a:a:a:a:a:a:a:b:c:d)$. Since $7a+b+c+d=0$ and
      $a+b+c+d=0$, we obtain $a=0=b+c+d$ and $\eta=(0:0:0:0:0:0:0:b:c:-b-c)$.
      Since $b=c=0$ is impossible, w.l.o.g.~we put $b=1$. 
      By $P(0)=P(1)=P(c)=P(-1-c)=0$ we have $2\lambda=1$, hence
      \[
       (\alpha:\beta)=(1:0),
      \]
      and no further conditions on $c\in\mathbb{C}$\,. Thus, we have found  
      $120=\frac{10\cdot9\cdot8}{3!}$ singular lines.
 \item[Case 12] We have $\eta=(1:1:1:1:1:-1:-1:-1:-1:-1)$ and, thus, $C_2=C_4=10$, 
      $C_3=0$, and
      \[
       P(X)=X^4-10\lambda X^2+10\lambda -1=\big(X^4-1\big)-10\lambda\big(X^2-1\big)\,.
      \]
      Hence, $P(\pm1)=0$ holds for all $\lambda$. This means that every single point in the 
      $\Sigma_{10}$-orbit of $\eta$ is a singularity of each hyperquintic 
      \mbox{$Q=Q_{(\alpha:\beta)}$} in the $\Sigma_{10}$-symmetric family in $\mathbb{P}^8$. 
      For this reason, from now on we will call these points 
      \emph{generic singularities} (cf.~\cite{Schm}).
      The length of the $\Sigma_{10}$-orbit of $\eta$ is $126=\frac{1}{2}\cdot{10\choose5}$.
 \item[Case 18] Due to $4a+4b+c+d=a+b+c+d=0$ we immediately obtain $b=-a$ and $d=-c$, hence  
      $\eta=(a:a:a:a:-a:-a:-a:-a:c:-c)$. Since $a=0$ leads us back to case 4, we put $a=1$ and 
      find $C_2=2c^2+8$, $C_3=0$, $C_4=2c^4+8$, and $P$ appropriate. Via $P(\pm1)=P(\pm c)=0$ 
      we obtain 
      \[
       0=\big( 2\lambda(c^2+4)-(c^2+1)\big)(c+1)(c-1)\,.
      \]
      With $c=\pm1$ we are back in case 12, 
      so we assume $c\not=\pm1$. Thus, 
      \[ 
       0=(2\lambda-1)c^2+(8\lambda-1)\,.
      \]
      This equation has no solution for $\lambda=\frac{1}{2}$, 
      but for $\lambda\not=\frac{1}{2}$ we have 
      \begin{align*}
       \hspace{4cm}&0=\beta c^2+(3\alpha+4\beta)\,.  & \hspace{3cm} (*)
      \end{align*}
      Thus, $\eta=(1:1:1:1:-1:-1:-1:-1:c:-c)$ is a singular point of $Q_{(\alpha:\beta)}$\,, 
      $(\alpha:\beta)\not=(1:0)$\,, for all $c\in\mathbb{C}$\, that satisfy $(*)$. 
      
      \noindent
      There are \mbox{$3150=\frac{1}{2}\cdot{10\choose4}{6\choose4}$} elements in the  
      $\Sigma_{10}$-orbit of $\eta$, 
      which we will also call 
      \emph{generic singularities} as well as in case 12  
      (cf.~\cite{Schm}).
      
      \noindent
      If $(\alpha:\beta)\in\{\,(5:-3),(4:-3)\,\}$\,, which means 
      $\lambda\in\{\,\frac{1}{5},\frac{1}{8}\,\}$, the solutions of $(*)$ are $c=\pm1$ and 
      $c=0$, respectively, so $\eta$ coincides with the singular points from case 12 or
      two orbit elements of $\eta$ merge to one of the singularities from case 17. 
      Hence, we have singularities that are worse than ordinary nodes. 
      A proof of this is given in section \ref{sec:ordinaryNodes}. 
      For this reason, we from now on will refer to 
      $(\alpha:\beta)\in\{\,(1:0),(5:-3),(4:-3)\,\}$\, or 
      $\lambda\in\{\,\frac{1}{2},\frac{1}{5},\frac{1}{8}\,\}$ as \emph{exceptional values}. 
      Case 7, however, already showed that $Q_{(1:0)}$ contains 120 singular lines.
\end{description}

We list the results of our investigation below. 
Table \ref{P8Tab1} shows the \emph{generic sin\-gu\-la\-ri\-ties}, 
which are contained in each hyperquintic of the family. For the \emph{exceptional values} $(\alpha:\beta)\in\left\{(1:0),(2:-1),(3:-2),(4:-3),(5:-3)\right\}$, the cor\-res\-ponding hyperquintics have singularities worse than ordinary nodes. 

In table \ref{P8Tab2} we list the parameter values, for which we have \emph{additional} orbits 
of singular points. Using computer algebra we can verify that all the \emph{additional} orbits consist only of ordinary nodes, if not stated otherwise. 

\begin{table}[bh]
 \begin{center}
 \small
  \begin{tabular}{|c|c|c|}
   \hline
   \rule{0pt}{2.3ex}orbit length & orbit element
   & case    \\ \hline \hline 
   \rule{0pt}{2.5ex}126          & $(1:1:1:1:1:-1:-1:-1:-1:-1)$
   & 12      \\ \hline 
                & \rule{0pt}{2.5ex}$(1:1:1:1:-1:-1:-1:-1:c:-c)$\,,
                &         \\ 
   \rb{3150}    & $\beta c^2+(3\alpha+4\beta)=0$           & \rb{18} \\ \hline
                & \rule{0pt}{2.5ex}$(1:1:1:-1:-1:-1:c:c:-c:-c)$\,,
                &         \\ 
   \rb{12600}   & $(\alpha+2\beta)c^2+(2\alpha+3\beta)=0$  & \rb{23} \\ \hline
  \end{tabular}
 \vspace{-0.2cm}
  \parbox{0.9\textwidth}{
   \caption{\emph{Generic singularities} in $\mathbb{P}^8$. Each hyperquintic 
            $Q_{(\alpha:\beta)}$ of the 1-parameter-family in $\mathbb{P}^8$ with 
            $(\alpha:\beta)$ not an \emph{exceptional value} contains these singular points. }
   \label{P8Tab1}
  }
 \end{center}
\end{table}

\begin{table}[hptb]
 \begin{center}
  \small
  \begin{tabular}{|c|c|c|}
   \hline
   $(\alpha:\beta)$ & orbit length          & orbit element\\ \hline \hline
                    &                       & $(1:1:1:1:1:b_1:b_1:b_2:b_2:3)$\,,\\
   \rb{$(3:-1)$}    & \rb{7560}             & $b_{1,2}=-2\pm\sqrt{-3}$\\ \hline
                    &                       & $(3:3:3:3:b_1:b_1:b_1:b_2:b_2:b_2)$\,,\\
   \rb{$(7:-5)$}    & \rb{4200}             & $b_{1,2}=-2\pm\sqrt{-3}$\\ \hline
                    &                       & $(1:1:1:1:1:1:-5:-5:c_1:c_2)$\,,\\
   \rb{$(51:-25)$}  & \rb{2520}             & $c_{1,2}=2\pm\sqrt{-7/5}$\\ \hline
                    &            
         & $(2\!:\!2\!:\!2\!:\!2\!:\!2\!:\!2b\!:\!2b\!:\!2b\!:\!-5\!-\!3b\!:\!-5\!-\!3b)$\\
   \rb{$(280\!:\!-163\pm\!3\sqrt{65})$}
                    & \rb{2520}             & $b=\frac{3\pm\sqrt{65}}{2}$\\ \hline
                    &                       & $(1:1:1:1:1:b:b:b:b:-4b-5)$\,,\\
   \rb{$(30:-13\mp\sqrt{85})$}
                    & \rb{1260}             & $b=\frac{-29\pm\sqrt{85}}{14}$\\ \hline
                    &                       & $(1:1:1:1:1:1:b_1:b_1:b_2:b_2)$\,,\\
   \rb{$(4:-1)$}    & \rb{1260}             & $b_{1,2}=\frac{-3\pm\sqrt{-7}}{2}$\\ \hline
                    &                       & $(1:1:1:1:1:1:b:b:b:-3b-6)$\,,\\
   \!\rb{$(21\!:\!2b^3\!+\!25b^2\!+\!86b\!+\!76)$}\!
                    & \rb{840}              & $2b^4+25b^3+93b^2+139b+77=0$\\ \hline 
                    &                       & $(1:1:1:1:1:1:1:b:b:-2b-7)$\,,\\
   \!\rb{$(84\!:\!175\!-\!3b(b^2\!+\!5b\!-\!13))$}\!
                    & \rb{360}              & $3b^4+39b^3+189b^2+413b+364=0$\\ \hline 
   $(50:-37)$       & 210                   & $(2:2:2:2:2:2:-3:-3:-3:-3)$\\ \hline
   $(175:-117)$     & 120                   & $(3:3:3:3:3:3:3:-7:-7:-7)$\\ \hline
                    &                       & $(1:1:1:1:1:1:1:1:b_1:b_2)$\,,\\
   \rb{$(3:7)$}     & \rb{90}               & $b_{1,2}=-4\pm\sqrt{-11}$\\ \hline
   $(100:-49)$      & 45                    & $(1:1:1:1:1:1:1:1:-4:-4)$\\ \hline
   $(75:7)$         & 10                    & $(1:1:1:1:1:1:1:1:1:-9)$\\ \hline
                    &                       & $(0:0:0:0:0:0:0:b:c:d)$\,,\\
   \rb{$(1:0)$}     & \rb{120 lines}        & $b+c+d=0$\\ \hline
                    &                       & $(0:0:0:0:b:b:c:c:d:d)$\,,\\
   \rb{$(2:-1)$}    & \rb{3150 lines}       & $b+c+d=0$\\ \hline
                    &                       & $(a:a:a:b:b:b:c:c:c:0)$\,,\\
   \rb{$(3:-2)$}    & \rb{2800 lines}       & $a+b+c=0$\\ \hline
                    & 1575 ($D_4$)          & \\
   \rb{$(4:-3)$}    & (Remark \ref{P8D4}) 
                                            & \rb{$(1:1:1:1:-1:-1:-1:-1:0:0)$}\\ \hline
                    & 126                   & \\
   \rb{$(5:-3)$}    & (Remark \ref{DelPezzo})
                                            & \rb{$(1:1:1:1:1:-1:-1:-1:-1:-1)$}\\ \hline
                    & hypersurface          & \\
   \rb{$(0:1)$}     & $S_2=S_3=0$           & \\ \hline
  \end{tabular}
  \vspace{-0.2cm}
  \parbox{0.9\textwidth}{
   \caption{Parameter values, for which we have \emph{additional} orbits of singular points. 
            By using computer algebra we can verify that only ordinary nodes are contained in 
            these orbits, if not stated otherwise. }
   \label{P8Tab2}
  }
 \end{center}
\end{table}

As we will see in the next section, all the \emph{generic singularities} are ordinary nodes. 
Moreover, for $(\alpha:\beta)=(3:-1)$, which corresponds to the longest orbit of \emph{additional} singular points, we find the best hyperquintic in the $\Sigma_{10}$-symmetric family in $\mathbb{P}^8$.

\begin{samepage}
\begin{satz}
\label{Rekord}
The hyperquintic $Q_{(3:-1)}$, given by
\[
 3\cdot S_5 +(-1)\cdot S_2 S_3\,=\,S_1\,=\,0,
\]
where $S_i,\; i=1,2,3,5,$ is the $i$-th elementary-symmetric polynomial in 10 variables, has exactly 23436 ordinary nodes and no further singularities. 
\end{satz}
\vspace{-8mm}
\begin{flushright}$\square$\end{flushright}
\end{samepage}

\section{Ordinary Nodes}
\label{sec:ordinaryNodes}
To show that all the isolated singularities are ordinary nodes, we use
the Hessian criterion, i.e.~we show $\det(\Hess_f(y))\not=0$, where
$\Hess_f=\left(\frac{\partial^2f}{\partial x_i \partial
    x_j}\right)_{\!ij}$ is the Hessian of $f$, $f=0$ is the affine
equation of the hyperquintic $Q_{(\alpha:\beta)}$ in an appropriate
affine chart, and $y$ is the singular point in this chart. 

\vspace{3mm}

Modulo $S_1$ one has 
\begin{align*}
 F:= & \,\textstyle 
    F_{(1:\beta)}=S_5+\beta\cdot S_2 S_3 
    =
    \frac{1}{5}C_5-\frac{1+\beta}{6}C_2 C_3\\
   = & \,\frac{1}{5}\bigg(\sum_{i=0}^{8} x_i^5-g(x)^5 \bigg)
         -\frac{1+\beta}{6}\bigg(\sum_{i=0}^{8} x_i^2+g(x)^2\bigg)
                           \bigg(\sum_{i=0}^{8} x_i^3-g(x)^3\bigg),
\end{align*}
where $g(x):=x_0+\hdots+x_8$.
We consider the isolated singularities in affine charts 
$\mathbb{A}^{8}_{\,i}$, $i\in\{\,0,\hdots,8\,\}$, given by 
\begin{align*}
 \mathbb{A}^{8}_{\,i}:=\big\{(x_0,\hdots,x_{i-1},x_{i+1},\hdots,x_8)|
      & (x_0:\hdots:x_{i-1}:1:x_{i+1}:\hdots:x_9)\in\mathbb{P}^8,\\
      & \hspace{1mm}\, x_9=-(x_0+\hdots+x_8)\,\big\}.
\end{align*}
Those charts cover the projective space $\mathbb{P}^{8}$, so that we find all the isolated singularities in at least one chart $\mathbb{A}^{8}_{\,i}$. 
In our case it is even sufficient to check only one chart, w.l.o.g.~$\mathbb{A}^{8}:=\mathbb{A}^{8}_{\,0}$, since no coordinate of 
our isolated singularities is zero.
Defining $h(x):=1+x_1+\hdots+x_{8}$, we obtain 
{\small
\begin{align*}
 f:=&f(x_1,\hdots,x_{8}):=\,F\big(1,x_1,\hdots,x_{8},-(1+x_1+\hdots+x_{8})\big)\\
   =&\,\frac{1}{5}\bigg(1+\sum_{k=1}^{8}x_k^5-h(x)^5\bigg)
     -\frac{1+\beta}{6}\bigg(1+\sum_{k=1}^{8}x_k^2+h(x)^2\bigg)
                       \bigg(1+\sum_{k=1}^{8}x_k^3-h(x)^3\bigg).
\end{align*}
}
Thus, it holds for the partial derivatives $f_i=\frac{\partial f}{\partial x_i}$, $i=1,\hdots,8$, of $f$ 
{\small
\begin{align*}
 f_i&=x_i^4-h(x)^4
      {\textstyle -\frac{1+\beta}{6}\bigg[2(x_i+h(x))+3(x_i^2-h(x)^2)}
           +2\cdot\sum_{k=1}^{8}x_k^3(x_i+h(x))\\
    &\hspace{3cm}+3\cdot\sum_{k=1}^{8}x_k^2(x_i^2-h(x)^2)-2x_ih(x)^3+3x_i^2h(x)^2-5h(x)^4\bigg]
\end{align*}
}

and for the second partial derivatives $f_{ii}=\frac{\partial^2 f}{\partial x_i^2}$ and $f_{ij}=\frac{\partial^2 f}{\partial x_i \partial x_j}$, $i\not=j$, 
\begin{align*}
 f_{ii}&=4x_i^3-4h(x)^3{\textstyle -\frac{1+\beta}{6}}\bigg[4+6x_i-6h(x)+12x_i^3+12x_i^2h(x)
           +4\cdot\sum_{k=1}^{8}x_k^3\\
       &\hspace{4.4cm}-6x_ih(x)^2+6(x_i-h(x))\sum_{k=1}^{8}x_k^2-22h(x)^3\bigg],\\
 f_{ij}&=-4h(x)^3\textstyle -\frac{1+\beta}{6}\bigg[2-6h(x)+6x_ix_j(x_i+x_j)
                                  +6h(x)(x_i^2+x_j^2)\\
       &\hspace{2.2cm}-6h(x)^2(x_i+x_j)+2\cdot\sum_{k=1}^{8}x_k^3
                                  -6h(x)\cdot\sum_{k=1}^{8}x_k^2-20h(x)^3\bigg].
\end{align*}

In the following subsections, we first check that all \emph{generic singularities} are ordinary nodes. Then we verify that the longest orbit of length 7560 of the \emph{additional} singularities of $Q_{(3:-1)}$ consists only of ordinary nodes. 

\subsection{The 126 generic Nodes}
We consider $\eta:=(1:1:1:1:1:-1:-1:-1:-1:-1)$ with its 126 orbit
elements; due to our choice of the affine chart $\mathbb{A}^{8}$ and
$S_1=0$, we evaluate the Hessian $\Hess_f=\big(f_{ij}\big)$ in
$y:=(1,1,1,1,-1,-1,-1,-1)$.  
With $h(y)=1$, we obtain 
\begin{align*}
 f_{ii}(y)=\left\{
            \begin{array}{r@{,\;}l}
                           0\;  & \text{ if }i\leq 4,\\
                   12 + 20\beta & \text{ if }i>4, 
            \end{array}
           \right.
& \quad\text{ and }\quad &
 f_{ij}(y)=6 + 10\beta \text{ for all }i\not=j.
\end{align*}

Thus,
\[
 \Hess_f(y)=\big(6 + 10\beta\big)\cdot
            \left(
             \begin{array}{cccccccc}
              0      & 1      & \cdots & \cdots & \cdots & \cdots & \cdots & 1\\
              1      & 0      & \ddots &        &        &        &        & \vdots\\
              \vdots & \ddots & 0      & \ddots &        & 1      &        & \vdots\\
              \vdots &        & \ddots & 0      & \ddots &        &        & \vdots\\
              \vdots &        &        &\ddots  & 2      & \ddots &        & \vdots\\
              \vdots &        & 1      &        & \ddots & 2      & \ddots & \vdots\\
              \vdots &        &        &        &        & \ddots & 2      & 1\\
              1      & \cdots & \cdots & \cdots & \cdots & \cdots & 1      & 2\\
             \end{array}
            \right).
\]
The determinant of the righthand matrix is 1, hence 
\[
 \det(\Hess_f(y))\not=0 \;\text{ for all }\textstyle
 \beta\not=-\frac{3}{5}. 
\]
But $(\alpha:\beta)=(1:-\frac{3}{5})=(5:-3)$ is one of the \emph{exceptional values}, hence all of the 126 orbit elements of $\eta$ are ordinary nodes of $Q_{(\alpha:\beta)}$, $(\alpha:\beta)\not=(5:-3)$, $\alpha\not=0$. 
For $(\alpha:\beta)=(5:-3)$ we have singularities worse than 
\label{NodesBeweis}
ordinary nodes. 

\subsection{The 3150 generic Nodes}
Now consider $\eta:=(1:1:1:1:-1:-1:-1:c:-c:-1)$ and its orbit elements, 
where $\beta c^2+(3+4\beta)=0$, in the affine chart $\mathbb{A}^{8}$. We put 
$ y:=(1,1,1,-1,-1,-1,c,-c)$ and obtain
\begin{align*}
 h(y)=(1+y_1+\hdots+y_{8})=1, \quad & \quad
 \sum_{k=1}^{8}y_k^3=0,       \quad &
 \sum_{k=1}^{8}y_k^2=6+2c^2. 
\end{align*}
Let 
\begin{align*}
 b_1:= & \, -4+(1+\beta)(8+2c^2)\,,          \\
 b_4:= & \; 6+10\beta\,,                     \\
 b_5:= & \; 4c^3-4-(1+\beta)(4c^3+6c-10)\,,  \\
 b_6:= & \, -4c^3-4+(1+\beta)(4c^3+6c+10)\,,
\end{align*}
then we have 
\begin{center}
 \begin{tabular}{ll}
  $f_{11}(y)=f_{22}(y)=f_{33}(y)=0,$ & $f_{44}(y)=f_{55}(y)=f_{66}(y)=2b_1,$ \\
  $f_{77}(y)=b_5,$                   & $f_{88}(y)=b_6$, 
 \end{tabular}
\end{center}
and for $i\not= j$ 
\[
 f_{ij}(y)=\left\{
            \begin{array}{r@{,\quad}l}
              b_1 & y_i=y_j=+1,\\
              b_1 & y_i=y_j=-1,\\
              b_1 & y_i=+1,\;y_j=-1,\\
              b_4 & y_i=+1,\;y_j=+a,\\
              b_4 & y_i=+1,\;y_j=-a,\\
              b_1 & y_i=-1,\;y_j=+a,\\
              b_1 & y_i=-1,\;y_j=-a,\\
              b_4 & y_i=+a,\;y_j=-a. 
            \end{array}
           \right.
\]

Hence, for the Hessian $\Hess_f(y)$ we have 
\[
 \Hess_f(y)=\left(
             \begin{array}{ccc|ccc|c|c}
              0   & b_1 & b_1 & b_1  & b_1  & b_1  & b_4 & b_4 \\
              b_1 & 0   & b_1 & b_1  & b_1  & b_1  & b_4 & b_4 \\
              b_1 & b_1 & 0   & b_1  & b_1  & b_1  & b_4 & b_4 \\ \hline
              b_1 & b_1 & b_1 & 2b_1 & b_1  & b_1  & b_1 & b_1 \\
              b_1 & b_1 & b_1 & b_1  & 2b_1 & b_1  & b_1 & b_1 \\
              b_1 & b_1 & b_1 & b_1  & b_1  & 2b_1 & b_1 & b_1 \\ \hline
              b_4 & b_4 & b_4 & b_1  & b_1  & b_1  & b_5 & b_4 \\ \hline
              b_4 & b_4 & b_4 & b_1  & b_1  & b_1  & b_4 & b_6 \\
             \end{array}
            \right).
\]
Performing row and column transformations, one easily finds 
\[
 \det(\Hess_f(y)) = 2^8\cdot c^2\cdot (c^2-1)^8\cdot\frac{3^2}{(c^2+4)^2}\,.
\]
The denominator is not zero, since this would lead to a contradiction with the constraint on $c$. 
So the determinant only vanishes for $c\in\{0,\pm1\}$. But $c$ takes these values only for 
\label{NodesBeweis2}
\mbox{$(\alpha:\beta)\in\{(5:-3),(4:-3)\}$}, which are \emph{exceptional values}.
Then we have singularities worse than ordinary nodes, due to certain merging singularities. 
For other values of $(\alpha:\beta)$, $\alpha\not=0$, all the 3150 orbit elements of $\eta$ are ordinary nodes. 

\subsection{The 12600 generic and the 7560 additional Nodes}
For the 12600 orbit elements of \mbox{$(1:1:1:-1:-1:-1:c:c:-c:-c)$} with $(1+2\beta)\cdot c^2+(2+3\beta)=0$ as well as for the 7560 \emph{additional} orbit elements of \mbox{$(1:1:1:1:1:b_1:b_1:b_2:b_2:3)$} with \mbox{$b_{1,2}=-2\pm\sqrt{-3}$,} 
the procedure is exactly the same. For the latter case, we take $\beta=-\frac{1}{3}$ into account, since it is an \emph{additional} orbit of singularities for $Q_{(3:-1)}=Q_{(1:-\frac{1}{3})}$. 
Thus, we find $\det(\Hess_f(y_1))\not=0$ 
for $(\alpha:\beta)\not=(5:-3),(3:-2)$, and
$\det(\Hess_f(y_2))\not=0$, where $y_1:=(1,1,-1,-1,c,c,-c,-c)$ and
$y_2:=(1,1,1,1,b_1,b_1,b_2,b_2)$.  
So the 12600 \emph{generic} and the 7560 \emph{additional} orbit elements of the corresponding singularities are ordinary nodes. The latter ones are contained only in $Q_{(3:-1)}$. 

\section{Concluding Remarks}
   \label{sec:conclRem}
We have proved in the previous sections that the hyperquintic
$Q_{(3:-1)}$ in $\mathbb{P}^8$ has 23436 ordinary nodes and no further
singularities. We now briefly discuss the cases $(\alpha:\beta)$,
where $Q_{(\alpha:\beta)}$ has higher singularities (see table
\ref{P8Tab2}).  
Moreover, we look at the generalization of our approach
to $\mathbb{P}^n$ and compare it to another construction of
hyperquintics in $\mathbb{P}^n$ with many nodes.  

\subsection{Some Cases with higher Singularities}
\label{sec:higherSings}
\begin{bemerkung}
\label{DelPezzo}
For $(\alpha:\beta)=(5:-3)$, 25 respectively 100 orbit elements of the
\emph{generic}  
singularities from cases 18 and 23, respectively, 
coincide with one appropriate orbit element of 
$
 (1:1:1:1:1:-1:-1:-1:-1:-1).
$
Thus, 126 singularities with Milnor number $256=2^8$ are created. The tangent cone of $Q_{(5:-3)}$ is a smooth cubic. 

We have similarities to this in the case of $\mathbb{P}^4$; here, $Q_{(3:-1)}$ has exactly the 10 orbit elements of $(1:1:1:-1:-1:-1)$ as singular points, they are called \emph{Del Pezzo Nodes} in \cite{vStr}. They have a Milnor number of \mbox{$16=2^4$,} the tangent cone of $Q_{(3:-1)}$ is a smooth cubic as well.
\end{bemerkung}

\begin{bemerkung}
\label{P8D4}
Besides the two orbits of \emph{generic} ordinary nodes, the hyperquintic $Q_{(4:-3)}$ in $\mathbb{P}^8$ has one more orbit with 1575 isolated singularities of type $D_4$, namely the orbit elements of $(1:1:1:1:-1:-1:-1:-1:0:0)$. 
See also \cite{Schm}.
\end{bemerkung}

\subsection{The Pentagon Construction}
\begin{figure}[t]
 \begin{center}
  \begin{minipage}[t]{10.3cm}
   \begin{tabular}{ccc}
    \includegraphics[width=4.5cm]{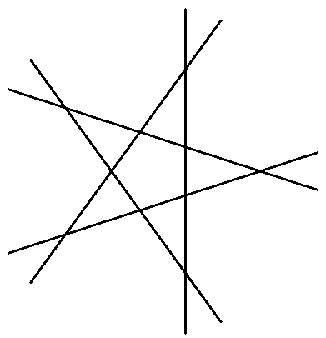} & \qquad &  \includegraphics[width=4.5cm]{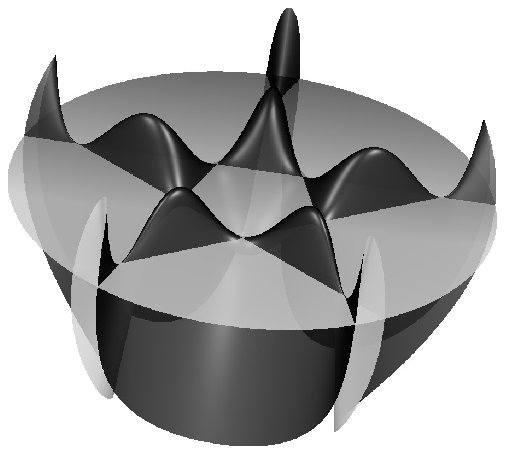}
   \end{tabular}
  \end{minipage}
  \caption{$R_5(x,y)=0$ in $\mathbb{P}^2$ and $z-R_5(x,y)=0$ in $\mathbb{P}^3$.}
  \label{D5eck}
 \end{center}
\end{figure}

Theorem \ref{Rekord} improves the previously best known lower bound of
23126 for the maximum number of ordinary nodes a hyperquintic in
$\mathbb{P}^8$ can have. 
The hypersurface corresponding to that previous lower bound was obtained
with an approach based on a generalization 
of constructions by Givental for cubic hypersurfaces and Hirzebruch
for quintics in $\mathbb{P}^4$ (cf.~\cite{AGZV}, \cite{Hir}, and
\cite{Lab}, sections 3.8--3.12).  
The basic idea is the usage of several polynomials of degree $5$ in two
variables, which have only a small number of critical values, to
construct hyperquintics with many nodes. More precisely, one considers
regular pentagons  
\[
 R_5(x,y):=x^5-10x^3y^2+5xy^4-5x^4-10x^2y^2-5y^4+20x^2+20y^2-16
\] 
in the plane. 
These can be normalized such that their critical values are $0$ and
$\pm1$ (see figure \ref{D5eck}).  
Then, Givental's equations for cubics can be transferred word by word
to obtain an affine equation for hyperquintics in $\mathbb{P}^n$ with
many singular points (all are nodes): 
\[
\sum_{j=0}^{\lfloor\frac{n}{2}\rfloor-1} 
(-1)^{j\cdot(1+(n\bmod2))}\widetilde{R}_5(x_{2j},x_{2j+1})= 
-(n\bmod2)\frac{T_5(x_{n-1})-1}{2}\,, 
\]
where $T_5(z):=16z^5-20z^3+5z$ denotes the Tchebychev polynomial of
degree $5$ with two critical values $\pm1$ and $\widetilde{R}_5(x,y)$ is
the normalized pentagon with critical value $+1$ over the origin.  
For a comparison of the resulting hyperquintics obtained by this
method to our $\Sigma_{n+2}$-symmetric approach see table
\ref{FuenfeckTab1}. 

\subsection{The $\Sigma_{n+2}$-symmetric Approach}
We performed further experiments for some $n\not=4$, and it seems to
us that the $\Sigma_{n+2}$-symmetric construction yields fewer nodes
than the pentagon construction.  
\begin{table}[bhtp]
 \begin{center}
 \small
 \begin{tabular}{|r|r|r|r|}
  \hline
           & \multicolumn{2}{c|}{\rule{0pt}{2.2ex}number of ordinary nodes} &    \\ \cline{2-3}
  \rb{$n$} &
  \multicolumn{1}{c|}{\rule{0pt}{2.3ex}$\Sigma_{n+2}$-symmetric
    approach}  
           & \multicolumn{1}{c|}{pentagon construction} &
           \multicolumn{1}{c|}{\rb{$\Ar_n(5)$}}\\ 
  \hline \hline
  \rule{0pt}{2.5ex}3   &            20   &        31   &        31  \\
  \hline 
  \rule{0pt}{2.5ex}4   &      {\bf 130}  &       126   &       135  \\
  \hline 
  \rule{0pt}{2.5ex}5   &           210   &       420   &       456  \\
  \hline 
  \rule{0pt}{2.5ex}6   &          1505   &      1620   &      1918  \\
  \hline 
  \rule{0pt}{2.5ex}8   &    {\bf 23436}  &     23126   &     27876  \\
  \hline  
  \rule{0pt}{2.5ex}10  &        296604   &    325580   &    411334  \\
  \hline 
 \end{tabular}
 \vspace{0.2cm}
 \parbox{0.9\textwidth}{
  \caption{Comparison of our $\Sigma_{n+2}$-symmetric approach 
           and the pentagon construction in $\mathbb{P}^n$ for some
           $n$. 
          } 
  \label{FuenfeckTab1}
 }
 \end{center}
\end{table}
Indeed, the best hyperquintic $Q_{(7:-4)}$ in $\mathbb{P}^5$ contains
only 210 ordinary nodes (cf.~table \ref{FuenfeckTab1}); 
in $\mathbb{P}^3$, the best hyperquintic $Q_{(2:1)}$ has only 
20 ordinary nodes. 
For $n=6$ and $n=10$ we obtained 1505 respectively 296604 ordinary
nodes for the best examples. 
In $\mathbb{P}^4$ and $\mathbb{P}^8$, the $\Sigma_{n+2}$-symmetric approach 
yields hyperquintics with a higher number of ordinary nodes than the 
pentagon construction (cf.~table \ref{FuenfeckTab1}). 

We did not look at other $n$ in detail, we only verified that the number of \emph{generic} nodes of the $\Sigma_{n+2}$-symmetric approach is less than the number of nodes obtained by using the pentagon construction. It is possible that for certain $n$ and $(\alpha:\beta)$ the $\Sigma_{n+2}$-symmetric construction is better.

\begin{appendix}
\section*{Appendix}
\markboth{\sl APPENDIX}{\sl APPENDIX}
\thispagestyle{myheadings}

\subsection*{The Case Analysis}
   \label{case_analysis}
Here we list the remaining cases of the case analysis in section \ref{family_hyperquintics}.
\begin{description}
 \item[Case 3] Assume that $\eta=(1:1:1:1:1:1:1:1:-4:-4)$. Hence $C_2=40$, $C_3=-120$, 
      $C_4=520$, and
      \[
       P(X)=\textstyle X^4-40\lambda X^2+80\lambda X+160\lambda -52\,.
      \]
      Via $P(1)=P(-4)=0$ we get $\lambda=\frac{51}{200}$, thus
      \[
       (\alpha:\beta)=(100:-49)\,.
      \]
      The length of the $\Sigma_{10}$-orbit of $\eta$ is $45={10\choose2}$\,.
 \item[Case 4] Consider $\eta=(a:a:a:a:a:a:a:a:b:-8a-b)$. For $a=0$ we get 
      $\eta=(0:0:0:0:0:0:0:0:1:-1)$, hence $C_2=C_4=2$, $C_3=0$, and
      \[
       \textstyle P(X)=X^4-2\lambda X^2+\frac{1}{5}(2\lambda-1)\,.
      \]
      Requiring $P(0)=P(\pm 1)=0$ leads to $\lambda=\frac{1}{2}$, hence
      \[
       (\alpha:\beta)=(1:0)\,.
      \]
      The length of the $\Sigma_{10}$-orbit of $\eta$ is $45=\frac{1}{2}\cdot10\cdot9$.
      
      \noindent
      For $a\not=0$ we have $\eta=(1:1:1:1:1:1:1:1:b:-8-b)$, hence
      $C_2=2b^2+16b+72$, $C_3=-24b^2-192b-504$, $C_4=2b^4+32b^3+384b^2+2048b+4104$, and $P$ 
      appropriate. By requiring $P(1)=P(b)=P(-8-b)=0$, we obtain three equations for 
      $\lambda$: 
      \begin{align*}
       0=&\lambda\cdot\big(\underbrace{2b^4+32b^3+342b^2+1712b+3912}_{=:a_{11}} \big)\\
         &\qquad +1\cdot\big(\underbrace{-b^4-16b^3-192b^2-1024b-2047}_{=:a_{12}} \big)\\
       0=&\lambda\cdot\big(\underbrace{-8b^4+32b^3+552b^2+2832b+2592}_{=:a_{21}} \big)\\
         &\qquad +1\cdot\big(\underbrace{4b^4-16b^3-192b^2-1024b-2052}_{=:a_{22}} \big)\\
       0=&\lambda\cdot\big(\underbrace{-8b^4-288b^3-3288b^2-16528b-33888}_{=:a_{31}} \big)\\
         &\qquad +1\cdot\big(\underbrace{4b^4+144b^3+1728b^2+9216b+18428}_{=:a_{32}} \big)
      \end{align*}
      To have a unique solution for $\lambda$, the matrix 
      \[
       \left(
        \begin{array}{cc}
          a_{11} & a_{12}\\ a_{21} & a_{22}\\ a_{31} & a_{32}
        \end{array}
       \right)
      \]
      must have rank 1. 
      Thus, the three $2\times2$ minors
      \[
       \left|
        \begin{array}{cc}
           a_{11} & a_{12}\\ a_{21} & a_{22}
        \end{array}
       \right| \text{, }
       \left|
        \begin{array}{cc}
           a_{11} & a_{12}\\ a_{31} & a_{32}
        \end{array}
       \right| \text{, }
       \left|
        \begin{array}{cc}
           a_{21} & a_{22}\\ a_{31} & a_{32}
        \end{array}
       \right| 
      \]
      must vanish, which leads to
      \[
       0=(b+9)(b-1)(b+4)(b^2+8b+27)\,.
      \]
      The solutions 1 and $-9$ lead us back to case 2, $b=-4$ to case 3. 
      If we take one of the two roots $b=-4\pm\sqrt{-11}$ of the remainig factor, 
      $-8-b$ is the other one. Thus, the length of the $\Sigma_{10}$-orbit of  
      $\eta$ is $90=10\cdot9$. 
      
      \noindent
      Such an $\eta$ yields $C_2=18$, $C_3=144$, $C_4=-1350$, and
      \[
       P(X)=X^4-18\lambda X^2-96\lambda X+\textstyle\frac{162}{5}\lambda+135\,.
      \]
      Requiring $P(1)=P(-4\pm i\sqrt{11})=0$ leads to $3\lambda=5$ and 
      \[
       (\alpha:\beta)=(3:7)\,.
      \]
 \item[Case 5] Assume that $\eta=(3:3:3:3:3:3:3:-7:-7:-7)$. Hence
      $C_2=210$, $C_3=-840$, $C_4=7770$, and 
      \[
       P(X)=X^4-210\lambda X^2+560\lambda X+4410\lambda-777\,.
      \]
      From $P(3)=P(-7)=0$ we get $\lambda=\frac{29}{175}$, hence
      \[
       (\alpha:\beta)=(175:-117)\,.
      \]
      The length of the $\Sigma_{10}$-orbit of $\eta$ is $120={10\choose3}$\,.
 \item[Case 6] Consider $\eta=(a:a:a:a:a:a:a:b:b:-7a-2b)$. For $a=0$ we put $b=1$, hence
      $C_2=6=-C_3$, $C_4=18$, and
      \[
       \textstyle P(X)=X^4-6\lambda X^2+4\lambda X+\frac{18}{5}\lambda-\frac{9}{5}\,.
      \]
      $P(0)=P(1)=P(-2)=0$ yields $\lambda=\frac{1}{2}$, 
      hence
      \[
       (\alpha:\beta)=(1:0)\,.
      \]
      The length of the $\Sigma_{10}$-orbit of $\eta$ is $360=10\cdot{9\choose2}$\,.

      \noindent
      For $a\not=0$, we have $\eta=(1:1:1:1:1:1:1:b:b:-7-2b)$. Thus, 
      \begin{align*}
       C_2&=6b^2+28b+56,\\
       C_3&=-6b^3-84b^2-294b-336,\\
       C_4&=18b^4+224b^3+1176b^2+2744b+2408,
      \end{align*}
      and $P$ appropriate.
      By requiring $P(1)=P(b)=P(-2b-7)=0$, we again obtain three equations for $\lambda$
      (cf.~case 4), hence a $3\times 2$-matrix, which must have rank 1 to have a unique 
      solution for $\lambda$. Thus, its three $2\times 2$-minors  must vanish and we find 
      \[
       0=(b-1)(3b+7)(b+4)(3b^4+39b^3+189b^2+413b+364)\,.
      \]
      The solutions $b\in\{\,1,-4,-\frac{7}{3}\,\}$ lead us back into cases 2, 3, and 5, 
      respectively. For $b$ a root of the remaining factor 
      and by $P(1)=P(b)=P(-2b-7)=0$, we obtain
      $\lambda=\frac{-1}{168}(3b^3+15b^2-39b-259)$, hence
      \[
       (\alpha:\beta)=(84:-3b^3-15b^2+39b+175)\,.
      \]
      The length of the $\Sigma_{10}$-orbit of $\eta$ is $360=10\cdot{9\choose2}$.
 \item[Case 8] We consider $\eta=(2:2:2:2:2:2:-3:-3:-3:-3)$ and obtain $C_2=-C_3=60$,
       $C_4=420$, and, thus, 
       \[
        P(X)=X^4-60\lambda X^2+40\lambda X+360\lambda -42\,.
       \]
       $P(2)=P(-3)=0$ leads to $\lambda=\frac{13}{100}$, hence  
       \[
        (\alpha:\beta)=(50:-37)\,.
       \]
       The length of the $\Sigma_{10}$-orbit of $\eta$ is $210={10\choose4}$.
 \item[Case 9] Assume that $\eta=(a:a:a:a:a:a:b:b:b:-6a-3b)$. For $a=0$, $b=1$ we find 
       $C_2=12$, $C_3=-24$, $C_4=84$, and
       \[
        P(X)=X^4-12\lambda X^2+16\lambda X+\textstyle\frac{72}{5}\lambda-\frac{42}{5}\,,
       \]
       but for no $\lambda$ does $P(0)=P(1)=P(-3)=0$ hold simultaneously.
       
       \noindent
       For $\eta=(1:1:1:1:1:1:b:b:b:-3b-6)$ we obtain 
       \begin{align*}
        C_2&=12b^2+36b+42,\\
        C_3&=-24b^3-162b^2-324b-210,\\
        C_4&=84b^4+648b^3+1944b^2+2592b+1302
       \end{align*}
       and $P$ appropriate.
       Equating $P(X)$ to zero for $X=1,b,-3b-6$ leads to an equation system for $\lambda$  
       again (cf.~cases 4 and 6), hence a $3\times2$-matrix, whose three $2\times2$-minors must 
       vanish to have a unique solution for $\lambda$. Thus, 
       \[
        0=(3b+7)(b-1)(2b+3)(b+2)(2b^4+25b^3+93b^2+139b+77)\,.
       \]
       The first three factors take us back into cases 5, 2, and 8, respectively, for $b=-2$ we 
       find $C_2=-C_3=18$, $C_4=54$, and
       \[
        \textstyle P(X)=X^4-18\lambda X^2+12\lambda X+\frac{162}{5}\lambda-\frac{27}{5}\,.
       \]
       Requiring $P(1)=P(-2)=P(0)=0$ yields $\lambda=\frac{1}{6}$, hence
       \[
       (\alpha:\beta)=(3:-2)\,.
       \]
      The length of the $\Sigma_{10}$-orbit of $\eta$ is $840=10\cdot{9\choose3}$.
      
      \noindent
      For $b$ a root of the remaining factor, we 
      find $\lambda=\frac{1}{42}(2b^3+25b^2+86b+97)$ by requiring $P(1)=P(b)=P(-3b-6)=0$. Thus, 
      \[
       (\alpha:\beta)=(21:2b^3+25b^2+86b+76)\,.
      \]
      The length of the $\Sigma_{10}$-orbit of $\eta$ also is $840=10\cdot{9\choose3}$ here.
 \item[Case 10] Assume that $\eta=(a:a:a:a:a:a:b:b:c:c)$ with $3a+b+c=0$\,.
      For $a=0$, $b=-c=1$ we find $C_2=C_4=4$, $C_3=0$, and
      \[
       P(X)=X^4-4\lambda X^2+\textstyle\frac{8}{5}\lambda-\frac{2}{5}\,.
      \]
      Via $P(0)=P(\pm1)=0$ we find $\lambda=\frac{1}{4}$, hence
      \[
       (\alpha:\beta)=(2:-1)\,.
      \]
      The length of the $\Sigma_{10}$-orbit of $\eta$ is
      $630=\frac{1}{2}\cdot{10\choose2}{8\choose2}$\,.
      
      \noindent
      For $a=1$, $c=-3-b$ we find $C_2=4(b^2+3b+6)$, $C_3=-6(3b^2+9b+8)$, 
      $C_4=4(b^4+6b^3+27b^2+54b+42)$, 
      and $P$ appropriate. The conditions on the coordinates of $\eta$ produce three equations 
      for $\lambda$. By the same method as in cases 4, 6, and 9, we obtain
      \[
       0=(b+4)(b-1)(2b+3)(b^2+3b+4)\,.
      \]
      The linear factors lead us back to cases 3 and 8. For a root 
      $b=\frac{-3\pm\sqrt{-7}}{2}$ of the last factor, $-b-3$ is the other one. Hence, the 
      length of the $\Sigma_{10}$-orbit of $\eta$ is $1260={10\choose2}{8\choose2}$.
      
      \noindent
      Such an $\eta=(2:2:2:2:2:2:2b:2b:-6-2b:-6-2b)$ 
      implies $C_2=32$, $C_3=192$, $C_4=-896$, and
      \[
       P(X)=X^4-32\lambda X^2-128\lambda X+\textstyle\frac{512}{5}\lambda+\frac{448}{5}\,.
      \]
      By $P(2)=P(2b)=P(-6-2b)=0$ we find $\lambda=\frac{3}{8}$, hence 
      \[
       (\alpha:\beta)=(4:-1)\,.
      \]
 \item[Case 11] Assume that $\eta=(a:a:a:a:a:a:b:b:c:d)$. 
      Due to $6a+2b+c+d=a+b+c+d=0$, we obtain $b=-5a$ and $d=4a-c$\,. 
      
      \noindent
      $a=0$ takes us back to case 4, so we put $a=1$. This leads to
      \begin{align*}
       C_2&=2c^2-8c+72,\\
       C_3&=12c^2-48c-180,\\
       C_4&=2c^4-16c^3+96c^2-256c+1512,
      \end{align*}
      and $P$ appropriate. By the conditions on $P$ we find 
      \begin{center}
       $5c^2-20c+27=0$ \quad and \quad $51\lambda =13$\,.
      \end{center}
      For a root $c=2\pm\sqrt{-\frac{7}{5}}$\,, the other one is $4-c$. Hence, we obtain 
      \[
       (\alpha:\beta)=(51:-25)
      \]
      and an orbit length of $\eta$ of $2520=10\cdot9\cdot{8\choose2}$.
 \item[Case 13] Consider $\eta=(a:a:a:a:a:b:b:b:b:-5a-4b)$. 
      For $a=0$ and $b=1$ 
      we find $C_2=20$, $C_3=-60$, $C_4=260$, and 
      \[
       P(X)=X^4-20\lambda X^2+40\lambda X+40\lambda -26\,.
      \]
      Requiring $P(0)=P(1)=P(-4)=0$ leads to a contradiction.
      
      \noindent
      For $a=1$ we have $C_2=20b^2+40b+30$, $C_3=-60b^3-240b^2-300b-120$, 
      $C_4=260b^4+1280b^3+2400b^2+2000b+630$, and $P$ appropriate.
      If we require $P(1)=P(b)=P(-4b-5)=0$, 
      we obtain the following equations:
      \begin{align*}
       0=&(b-1)(b+1)^3(2b+3)(7b^2+29b+27)\,,\\
       0=&(b+1)(450\lambda+182b^5+1055b^4+1795b^3+145b^2-2055b-1368)\,.
      \end{align*}
      The solutions 
      $b\in\{1,-1,-\frac{3}{2}\}$ of the first equation lead us back to 
      cases 2, 12, and 8, respectively, the roots $b=\frac{-29\pm\sqrt{85}}{14}$ of the 
      remaining factor together with the second equation imply 
      $\lambda=-\frac{7}{30}b-\frac{1}{5}=\frac{17\mp\sqrt{85}}{60}$, hence
      \[
       (\alpha:\beta)=(30:-13\mp\sqrt{85})\,.
      \]
      Here, the length of the $\Sigma_{10}$-orbit of $\eta$ is $1260=10\cdot{9\choose4}$\,.
 \item[Case 14] Assume that $\eta=(a:a:a:a:a:b:b:b:c:c)$. For $a=0$, $b=2$ we obtain 
       $C_2=-C_3=30$, $C_4=210$, and 
       \[
        P(X)=X^4-30\lambda X^2+20\lambda X+90\lambda-21\,.
       \]
       Requiring $P(0)=P(2)=P(-3)=0$, however, leads to a contradiction.
       
       \noindent
       For $\eta=(2:2:2:2:2:2b:2b:2b:-5-3b:-5-3b)$ we find $C_2=10(3b^2+6b+7)$, 
       $C_3=-30(b+7)(b+1)^2$, $C_4=10(21b^4+108b^3+270b^2+300b+133)$, and $P$ appropriate. 
       Requiring $P(2)=P(2b)=P(-5-3b)=0$ implies 
       \begin{align*}
        \hspace{0.3cm}&0=(b+1)(78400\lambda-186b^5+85b^4+4265b^3+6565b^2-6235b-24486)\,,\\
        \hspace{0.3cm}&0=(b-1)(b+1)^3(3b+7)(b^2-3b-14)\,.
       \end{align*}
       The roots $b\in\{\,1,-1,-\frac{7}{3}\,\}$ of the second equation take us back to 
       cases 3, 12, and 5, respectively. The roots $b=\frac{3\pm\sqrt{65}}{2}$ of the 
       remaining factor of the second equation, however, imply 
       $\lambda=\frac{117\pm3\sqrt{65}}{560}$ via the first equation, hence 
       \[
        (\alpha:\beta)=(280:-163\pm3\sqrt{65})\,.
       \]
       The length of the $\Sigma_{10}$-orbit of $\eta$ is 
       $2520={10\choose2}{8\choose3}$\,.
 \item[Case 15] Due to $5a+3b+c+d=a+b+c+d=0$ we may assume 
       $\eta=(a:a:a:a:a:-2a:-2a:-2a:c:a-c)$. For $a=0$ we are taken to case 4. 
       
       \noindent 
       For $a=1$ we find 
       $C_2=2(c^2-c+9)$, $C_3=3(c^2-c-6)$, $C_4=2(c^4-2c^3+3c^2-2c+27)$, 
       and $P$ appropriate. 
       But $P(1)=P(-2)=P(c)=P(1-c)=0$ imply $0=c(c-1)$ and $6\lambda = 1$\,,  
       and both $c=0$ and $c=1$ lead to case 9.
 \item[Case 16] Due to $5a+2b+2c+d=a+b+c+d=0$ we may assume 
       $\eta=(a:a:a:a:a:b:b:-4a-b:-4a-b:3a)$. Since $a=0$ leads to case 10, we put 
       $a=1$ and obtain $C_2=2(b^2+8b+23)$, $C_3=-24(b+2)^2$, 
       $C_4=4b^4+32b^3+192b^2+512b+598$, and $P$ appropriate. 
       Requiring $P(1)=P(b)=P(3)=P(-4-b)=0$, we find $b^2+4b+7=0$ and $\lambda=\frac{1}{3}$, 
       hence 
       \[
        (\alpha:\beta)=(3:-1)\,.
       \]
       For $b=-2\pm\sqrt{-3}$ one of the two solutions, $-4-b$ is the other one, so we find 
       $7560=10\cdot{9\choose2}{7\choose2}$ elements in the $\Sigma_{10}$-orbit of $\eta$.
 \item[Case 17] Assume that $\eta=(a:a:a:a:b:b:b:b:-2a-2b:-2a-2b)$. For $a=0$, $b=1$ we find 
      $C_2=-C_3=12$, $C_4=36$, and 
      \[
       P(X)=X^4-12\lambda X^2+8\lambda X+\textstyle\frac{72}{5}\lambda-\frac{18}{5}\,. 
      \]
      Via $P(0)=P(1)=P(-2)=0$ we immediately obtain $\lambda=\frac{1}{4}$, hence 
      \[
       (\alpha:\beta)=(2:-1)\,.
      \]
      The length of the $\Sigma_{10}$-orbit of $\eta$ is $3150={10\choose2}{8\choose4}$\,.
      
      \noindent
      For $a=1$ we find $C_2=4(3b^2+4b+3)$, $C_3=-12(b^3+4b^2+4b+1)$, 
      $C_4=4(9b^4+32b^3+48b^2+32b+9)$, and $P$ appropriate. 
      Requiring $P(1)=P(b)=P(-2b-2)=0$, we obtain 
      \begin{align*}
       \hspace{0.3cm}&0=b(b-1)(b+1)^3(3b+2)(2b+3)\,, \\
       \hspace{0.3cm}&0=400\lambda+210b^6+695b^5+556b^4-377b^3-742b^2-344b-100\,.  &(*)
      \end{align*}
      The case $b=0$ is checked above, whereas $b=1$ and $b\in\{-\frac{2}{3},-\frac{3}{2}\}$ 
      take us back to cases 3 and 8, respectively. 
      For $b=-1$ we obtain $\lambda=\frac{1}{8}$ via $(*)$, hence 
      \[
       (\alpha:\beta)=(4:-3)\,.
      \] 
      The length of the $\Sigma_{10}$-orbit of $\eta$ 
      is $1575=\frac{1}{2}\cdot{10\choose4}{6\choose4}$\,.
 \item[Case 19] We may assume $\eta=(3a:3a:3a:3a:b:b:b:-4a-b:-4a-b:-4a-b)$. For $a=0$, $b=1$ 
      we find $C_2=C_4=6$, $C_3=0$, and 
      \[
       P(X)=\textstyle X^4-6\lambda X^2+\frac{18}{5}\lambda-\frac{3}{5}\,.
      \]
      $P(0)=P(\pm1)=0$ leads to $\lambda=\frac{1}{6}$, hence 
      \[
       (\alpha:\beta)=(3:-2)\,.
      \] 
      The length of the $\Sigma_{10}$-orbit of $\eta$ is 
      $2100=\frac{1}{2}\cdot{10\choose3}{7\choose3}$.
      
      \noindent
      For $a=1$ we find 
       $C_2=6(b^2+4b+14)$, 
       $C_3=-12(3b^2+12b+7)$, 
       $C_4=6(b^4+8b^3+48b^2+128b+182)$,
      and $P$ appropriate. 
      Requiring $P(3)=P(b)=P(-4-b)=0$ leads to 
      \begin{align*}
       0&= (b-3)(b+2)(b+7)(b^2+4b+7)\,,\\
       0&= 4900\lambda-b^4-8b^3-6b^2+40b-581\,.
      \end{align*}
      The solutions $b\in\{\,-7,-2,3\,\}$ of the first equation take us back to cases 5 
      and 8, respectively. If $b$ is one of the two remaining roots $-2\pm\sqrt{-3}$ of the 
      first equation, $-4-b$ is the other one. Hence, the length of the $\Sigma_{10}$-orbit of 
      $\eta$ is $4200={10\choose3}{7\choose3}$\,. 
      Furthermore, by the second equation we obtain $\lambda=\frac{1}{7}$, hence 
      \[
       (\alpha:\beta)=(7:-5)\,.
      \]
 \item[Case 20] Due to $4a+3b+2c+d=a+b+c+d=0$ we have $c=-3a-2b$, $d=2a+b$, and  
      $\eta=(a:a:a:a:b:b:b:-3a-2b:-3a-2b:2a+b)$. Since $a=0$ takes us to case 17, we put $a=1$ 
      and obtain 
      \begin{align*} 
       C_2&=12b^2+28b+26,\\
       C_3&=-12b^3-66b^2-96b-42,\\
       C_4&=36b^4+200b^3+456b^2+464b+182,
      \end{align*} 
      and $P$ appropriate. Equating $P(X)$ to zero for $X=1,b,-3-2b,2+b$, we find 
      \begin{align*}
       0&= (b+1)^3(b+2)\,,\\
       0&= (b+1)(300\lambda+11b^3+18b^2-33b-100)\,.
      \end{align*}
      The two solutions $b\in\{-2,-1\}$ of the first equation lead us back to cases 9 and 
      12, respectively.
 \item[Case 21] Since we have $2a+b+c+d=a+b+c+d=0$, we immediately may assume 
      $\eta=(0:0:0:0:b:b:c:c:-b-c:-b-c)$. As $b$ and $c$ cannot vanish simultaneously, 
      w.l.o.g.~we put $b=1$ and find $C_2=4(c^2+c+1)$, $C_3=-6c(c+1)$, $C_4=4(c^2+c+1)^2$, and 
      $P$ appropriate. $P(0)=P(1)=P(c)=P(-1-c)=0$ imply $\lambda=\frac{1}{4}$, so 
      \[
       (\alpha:\beta)=(2:-1)\,.
      \] 
      We thus have found $3150=\frac{1}{3!}{10\choose2}{8\choose2}{6\choose2}$ singular lines,
      since there are no further conditions on $c\in\mathbb{C}$\,.
 \item[Case 22] Due to $3a+3b+3c+d=a+b+c+d=0$ we have $a+b+c=0=d$ and, hence, 
      $\eta=(a:a:a:b:b:b:-a-b:-a-b:-a-b:0)$. Since $a=0$ leads us back to case 19, we put $a=1$ 
      and find $C_2=6(b^2+b+1)$, $C_3=-9b(b+1)$, $C_4=6(b^2+b+1)^2$, and $P$ appropriate. 
      Requiring $P(0)=P(1)=P(b)=P(-1-b)=0$, we obtain $\lambda=\frac{1}{6}$, hence 
      \[
       (\alpha:\beta)=(3:-2)\,.
      \] 
      Since there are no further conditions on $b\in\mathbb{C}$\,, we again have found 
      $2800=\frac{1}{3!}{10\choose3}{7\choose3}{4\choose3}$ singular lines.
 \item[Case 23] Due to $3a+3b+2c+2d=a+b+c+d=0$ we have $b=-a$, $d=-c$ and, hence, may assume 
      $\eta=(a:a:a:-a:-a:-a:c:c:-c:-c)$. For $a=0$ we are back in case 10, so we put $a=1$. 
      Thus, we find $C_2=4c^2+6$, $C_3=0$, $C_4=4c^4+6$, and $P$ appropriate. Via 
      $P(\pm c)=P(\pm1)=0$ we obtain 
      \[
       0=\big(2\lambda(2c^2+3)-(c^2+1)\big)(c+1)(c-1)\,.
      \]
      With $c=\pm1$ we are back in case 12, so we assume $c\not=\pm1$. Thus, 
      \[
       0=(4\lambda-1)c^2+(6\lambda-1)\,.
      \]
      This equation has no solution for $\lambda=\frac{1}{4}$, but for 
      $\lambda\not=\frac{1}{4}$ we have 
      \begin{align*}
       \hspace{3.5cm}&0=(\alpha+2\beta)c^2+(2\alpha+3\beta)\,.  & \hspace{2.5cm} (*)
      \end{align*}
      So $\eta$ is a singular point of $Q_{(\alpha:\beta)}$, $(\alpha:\beta)\not=(2:-1)$, 
      for all $c\in\mathbb{C}$ that satisfy $(*)$ and we have found 
      \mbox{$12600=\frac{1}{2}\cdot{10\choose3}{7\choose3}{4\choose2}$} 
      more \emph{generic sin\-gu\-la\-ri\-ties} of $Q_{(\alpha:\beta)}$ in $\mathbb{P}^8$ 
      (cf.~\cite{Schm} and the cases 12 and 18). 
      
      \noindent
      If $(\alpha:\beta)\in\{\,(5:-3),(3:-2)\,\}$\,, which means 
      $\lambda\in\{\,\frac{1}{5},\frac{1}{6}\,\}$, the solutions of $(*)$ are $c=\pm1$ and 
      $c=0$, respectively, so two respective orbit elements of $\eta$ merge or they coincide 
      with the singular points from case 12. 
      Hence, we have singularities that are worse than ordinary nodes. 
      A proof of this is given in section \ref{sec:ordinaryNodes}. \\
      For this reason, we add $(\alpha:\beta)\in\{\,(2:-1),(5:-3),(3:-2)\,\}$\, and 
      $\lambda\in\{\,\frac{1}{4},\frac{1}{5},\frac{1}{6}\,\}$ to the \emph{exceptional values} 
      introduced in case 18. The cases 21 and 22, however, already showed that we have singular 
      lines contained in $Q_{(2:-1)}$ and $Q_{(3:-2)}$. 
\end{description}

\end{appendix}



\begin{thebibliography}{einrückung}
 \bibitem[AGZV]{AGZV} V.I.~Arnold, S.M.~Gusein-Zade, A.N.~Varchenko,
   {\it Singularities of differential maps,} Birkhäuser, 1985
 \bibitem[Bar]{Bar}   W.~Barth, {\it Two Projective Surfaces with Many
     Nodes, Admitting the Symmetry of the Icosahedron,} J.~Algebraic
     Geom.~{\bf 5} (1996), no.~1, 173--186  
 \bibitem[Gor]{Gor}   V.V.~Goryunov, {\it Symmetric Quartics with many Nodes,} Advances in 
                      Soviet Mathematics, Volume 21, 1994, 147--161
 \bibitem[Hir]{Hir}   F.~Hirzebruch, {\it Some Examples of Threefolds with trivial canonical 
                      bundle,} in Collected Works Vol. {\bf II} (1995), Springer Verlag, 
                      757--770
 \bibitem[JR]{JR}     D.B.~Jaffe, D.~Ruberman, {\it A Sextic Surface cannot have 66 Nodes,} 
                      J.~Algebraic Geom.~{\bf 6} (1997), no.~1, 151--168 
 \bibitem[Kal]{Kal}   T.~Kalker, {\it Cubic Fourfolds with 15 Ordinary Double Points,} 
                      Ph.~D.~thesis, Leiden, 1986
 \bibitem[Lab]{Lab}   O.~Labs, {\it Hypersurfaces with many Singularities,} Ph.~D.~thesis, 
                      Mainz, available from www.OliverLabs.net, 2005
 \bibitem[Schm]{Schm} O.~Schmidt, {\it Symmetrische Quintiken mit vielen Doppelpunkten,} 
                      Diploma thesis, Mainz, 2006
 \bibitem[vStr]{vStr} D.~van Straten, {\it A Quintic Hypersurface in $\mathbb{P}^4$ with 130
                      Nodes,} Topology {\bf 32} (1993), No. 4, 857--864
 \bibitem[Var]{Var}   A.N.~Varchenko, {\it On the Semicontinuity of the Spectrum and an Upper 
                      Bound for the Number of Singular Points of a Projective Hypersurface,}
                      J. Soviet Math. {\bf 270} (1983), 735--739
\end{thebibliography}
\end{document}